%% file: NewtonSebastianSubmit.tex
\documentclass[12pt]{amsart}

\usepackage{amsthm,amsfonts,amssymb,amsmath,color,graphicx,a4}

\author{Sebastian Mayer, Dierk Schleicher}
\title[Newton's Method for Entire Functions]
{Immediate and Virtual Basins of Newton's Method for Entire
Functions}
\address{{Lehrstuhl A f\"ur Mathematik, RWTH Aachen, D-52056
Aachen, Germany}}
\email{sebastian.mayer@mathA.rwth-aachen.de}
\address{{School of Engineering and Science, International University
Bremen, Postfach 750 561, D-28725 Bremen, Germany}}
\email{dierk@iu-bremen.de}



\date{\today}

\newtheorem{Theorem}{Theorem}[section]
\newtheorem{Lemma}[Theorem]{Lemma}

\newtheorem{Remark}[Theorem]{Remark}
\theoremstyle{remark}

\theoremstyle{definition}
\newtheorem{Definition}[Theorem]{Definition}

\newcommand{\ignore}[1]{}
\newcommand{\N}{\mathbb{N}}
\newcommand{\C}{\mathbb{C}}
\newcommand{\R}{\mathbb{R}}

\newcommand{\D}{\mathbb{D}}
\renewcommand{\S}{\mathbb{S}}
\newcommand{\Crit}{\textrm{Crit}}
\newcommand{\id}{\textrm{id}}

\renewcommand{\)}{\right)}

\newcommand{\image}[1]{\textrm{image}\left( #1 \right)}   

\newcommand{\ext}[1]{\textrm{ext}_{#1}}
\newcommand{\hide}[1]{}

\newcommand{\lineclear}
    {\rule{0pt}{0pt}\nopagebreak\par\nopagebreak\noindent}
\newcommand{\ovl}[1]{\overline{#1}}
\newcommand{\sm}{\setminus}
\newcommand{\remark}{\noindent{\sc Remark. }}

\begin{document}
\begin{abstract}
We investigate the well known Newton method to find roots of entire
holomorphic functions. Our main result is that the immediate basin of
attraction for every root is simply connected and unbounded. We also
introduce ``virtual immediate basins'' in which the dynamics converges
to infinity; we prove that these are simply connected as well.
\end{abstract}

\maketitle
\renewcommand{\thefootnote}{\fnsymbol{footnote}}
\renewcommand{\thefootnote}{\arabic{footnote}}

\section{Introduction}
\label{Introduction} 

Newton's method is one of the preferred methods to find roots of
differentiable maps: it often converges very fast and it is very
easy to implement. But there are problems: for example, even for
polynomials, there are open sets of initial conditions for which
the Newton map does not converge to any root.

There has been substantial progress understanding the
dynamics of Newton's method for finding roots of complex
polynomials: Przytycki~\cite{Feliks} has shown that all immediate
basins are simply connected and unbounded;
Shishikura~\cite{Mitsu} has shown more generally that if a
rational map has a multiply connected Fatou component, then it
must have two repelling or parabolic fixed points (which is
impossible for Newton maps of polynomials). Hubbard, Schleicher
and Sutherland~\cite{HSS} used this to find a rather small set of
starting points which together find all roots of a complex
polynomial; and in \cite{Dierk2} there is a (not very efficient)
bound on the number of iterations it takes to find all roots with
given accuracy. In a different spirit, Smale~\cite{Smale} has
shown that Newton's method is quite efficient from a
probabilistic point of view.

Newton's method for transcendental entire functions is less understood.
Bergweiler and Terglane~\cite{BT} have shown that Newton maps have no
multiply connected wandering domains, and in certain cases no wandering
domains at all. 

In this note, we extend Przytycki's result to entire holomorphic
functions: {\em for every root $\xi$ of a non-constant entire holomorphic
function $f$, the immediate basin of attraction for the Newton map
associated to $f$ is simply connected and unbounded}. This result goes
back to the Diploma thesis \cite{SebastianDiplom}. While the
result and part of its proof are in analogy to the polynomial
case, it turns out that for entire holomorphic functions, but not
for polynomials, ``virtual immediate basins'' play an important
role: these are domains in which the dynamics converges to
$\infty$ as if there was a root at $\infty$ (subject to further
conditions; see Definition~\ref{DefVirtualImmediateBasin}). We
show also that {\em every virtual immediate basin is simply
connected}. The combinatorial restrictions on ordinary and
virtual immediate basins are investigated in
\cite{JohannesNewton}.

It would be interesting to extend the ideas of \cite{Mitsu} to
the transcendental case, showing that all Fatou components are
simply connected; the case of wandering domains is treated in
\cite{BT}.

\bigskip
{\bf Acknowledgements.}  We would like to thank
 Johannes R\"uckert for his many helpful comments. We
 are grateful for the hospitality and the constructive atmosphere at the
Institut Henri Poincar\'e, Universit\'e Paris VI.
\section{Immediate Basins}
\label{Immediate basins}

Throughout this paper, let $f:\C \rightarrow \C$ be a nonlinear entire
holomorphic map and $N_f=\id - f/f'$ its associated Newton map. We will
be concerned with the set of points which converge to any given root $\xi$
of $f$. Clearly, the roots of $f$ are exactly the fixed points of
$N_f$ in $\C$, and these are attracting.

\begin{Definition}[Immediate basin]
\label{DefImmediateBasin}\lineclear
Let $\xi$ be an attracting fixed point of $N_f$. The \emph{basin of
attraction} of $\xi$ is the open set of all points $z$ such that
$(N_f^{\circ m}(z))_{m \in \N}$ converges to $\xi$. The connected
component containing $\xi$ of the basin is called the \emph{immediate
basin} of $\xi$.
\end{Definition}

Throughout this paper, we will fix a root $\xi$ of $f$ and denote its
immediate basin by $U$. In order to show that $U$ is both simply
connected and unbounded, we will construct a curve $\gamma\colon\R_0^+\to
U$ with $N_f(\gamma(t))=\gamma(t-1)$ for $t\geq 1$. If $U$ fails to be
full, then we can arrange things so that $\gamma(\R_0^+)$ is bounded
(and the same is obviously true if $U$ itself is bounded). In this case,
we show that $\lim_{t\to+\infty}\gamma(t)$ is a fixed point of $N_f$ in
$\partial U\cap \C$, and this will lead to a contradiction.

\subsection{An Exhaustion of Immediate Basins}
\label{Immediate basins:An exhaustion}

Let $\Crit(N_f)$ be the set of critical values of $N_f$ and 
\[
P_{U}:= {\bigcup_{m\ge 0} N_f^{\circ m}(\Crit(N_f)\cap U)}
\]
be the postcritical set restricted to critical values in $U$.
Since $\Crit(N_f)$ is countable, the set $P_U$ is countable as
well (but in general not closed). There is thus an open disk
$S_0=B_r(\xi)\subset U$ centered at $\xi$ such that
$\partial S_0\cap P_U=\emptyset$, and 
small enough such that $N_f(\overline{S_0}) \subset S_0$. For
every $k\in \N_0$ define $S_{k+1}$ to be the connected component
of $N_f^{-1}(S_k)$ containing $S_0$; then
$\ovl{S_{k+1}}$ is the connected component of $N_f^{-1}(\ovl{S_k})$
containing $S_0$.
\begin{Lemma}
\label{BasinIsCircleInverseImageIB}
The immediate basin satisfies
 $   U=\bigcup_{k\in \N} S_k\,\,.$
\end{Lemma}

\begin{proof}
Clearly, $U$ is open and $U':=\bigcup_{k\in \N} S_k$
is an open subset of $U$. Suppose there is a $z\in U \sm U'$.
Then there is an $M\in \N$ with
$N_f^{\circ M}(z)\in S_0$, so there is a connected neighborhood
$V\subset U$ of $z$ with $N_f^{\circ M}(V)\subset S_0$. For all $m\geq
M$, $z$ is by assumption in a component of $N_f^{-m}(S_0)$ different
from $S_m$, and so is $V$. Hence $U\setminus U'$ is open in
contradiction to the fact that $U$ is connected.
\end{proof}

Recall that a subset of $\C$ is called \emph{full} if its
complement has no bounded components. Clearly, if $U$ is not full, then
some of the $S_k$ are not full (if $U$ is not full, then it
contains a non-contractible loop which is compact and thus
contained in finitely many $S_k$).

\begin{Lemma}[Bounded $C_k$]
\label{LemBoundedSkIB}
If $S_M$ is not full but $S_{M-1}$ is, then all $S_m$ with $m<M$ are
bounded and homeomorphic to open disks.
\end{Lemma}

\begin{proof}
There is a bounded connected component $B$ of
$\C\sm S_M$. Its boundary $\partial B$ is a compact subset of
$\ovl{S_M}$, so $N_f(\partial B)$ is a compact subset of
$\ovl{S_{M-1}}$. There are no postcritical points in $\partial S_0$, so
$N_f^{\circ M}$ restricted to a neighborhood of $\partial B$ is a local
injection and $N_f^{\circ M}\colon\partial B\to\partial S_0$ is a
covering map. This implies that $\partial B$ and $N_f(\partial B)$ are
homeomorphic to circles. Since $N_f(\partial B)$ is a boundary component
of $S_{M-1}$ and $S_{M-1}$ is full, it follows that $S_{M-1}$ is
contained in the bounded complementary component of $N_f(\partial B)$,
so $S_{M-1}$ is bounded and homeomorphic to an open disk.

Clearly, all $S_m$ with $m<M$ are contained in $S_{M-1}$, hence
also bounded and simply connected.
\end{proof}
\subsection{Extending Paths Invariantly to Infinity}
\label{Immediate basins: Paths to infinity}

The goal of this section is the construction of curves
$\delta\colon\R_0^+\to U$ such that $N_f(\delta(t))=\delta(t-1)$
for $t\geq 1$.

\begin{Definition}[Extension of a curve]
If $\delta:[0,1]\rightarrow \C \setminus \Crit(N_f)$ is a curve
with $N_f(\delta(1))=\delta(0)$, then define its
\emph{(maximal) extension} $\ext{\delta}:[0,M_\delta) \rightarrow \C
\setminus \Crit (N_f)$ to be the unique curve with
\label{ext}
\begin{eqnarray*}
    \forall\, t \in [0,M_{\delta}-1): \quad N_f(\ext{\delta}(t+1))&=&\ext{\delta}(t)\\
    \ext{\delta} |_{[0,1]} \equiv \delta
\end{eqnarray*}
where $M_\delta$ is chosen maximal in $[1,\infty) \cup \lbrace
\infty \rbrace$.
\end{Definition}

\begin{Lemma}[Possibilities for the extension of a curve]
\label{LemExtCurvePossibilities}\lineclear
Given any curve $\delta\colon[0,1]\to\C\sm\Crit(N_f)$, then exactly one
of three cases occurs: (i) $M_\delta=\infty$, (ii) $M_{\delta}<\infty$ and
$\lim_{t \to M_{\delta}} \ext{\delta}(t)$ is a critical point of $N_f$,
or (iii) $M_{\delta}<\infty$, $\ext{\delta}(M_{\delta}-1)$ is an
asymptotic value, and ${\ext{\delta}(t)}\to\infty$ as $t\to
M_\delta$.
\end{Lemma}
\begin{proof}
Choose $T\in(0,M_\delta]$ with $T<\infty$ and suppose there
is a sequence $t_n\nearrow T$ with $\ext{\delta}(t_n)\to a\in\C$. If
$a$ is not a critical point, then there is a neighborhood $V$ of $a$
such that $N_f|_V$ is univalent, and it follows that $T<M_\delta$.
Therefore, if $M_\delta<\infty$, then $\ext{\delta}$ converges
either to a critical point or to infinity along an asymptotic
path.
\end{proof}

The following lemma is related to typical proofs of ``landing of
periodic dynamic rays'' for iterated polynomials. 

\begin{Lemma}[A homotopy class of unbounded curves]
\label{BoundedBackwardsImagesOfCurvesConverge}\lineclear
Let $W,W'\subset U\setminus S_0$ be two simply connected domains such
that $N_f(W') \subset W$ and $W' \subset W$. Let
\begin{equation*}
    S_{W',W}:=\{ \sigma : [0,1] \rightarrow W \textrm{ continuous,
} \sigma(1) \in
    W',\, \sigma(0)=N_f(\sigma(1))\} \,\, .
\end{equation*}
If there is an extension $\ext{\sigma}:[0,\infty)\rightarrow \C
\setminus \Crit(N_f)$ for all $\sigma\in S_{W,W'}$, then
every $\sigma\in S_{W,W'}$ has
$\lim_{t\to\infty}\ext{\sigma}(t)=\infty$.
\end{Lemma}

\begin{proof}
Set $Y:=\bigcup_{\sigma\in S_{W,W'}}\image{\ext{\sigma}}$. We use
$S_{W',W}$ to construct a sequence of local inverse mappings of
$N_f^{\circ m}$ on $W$. Note that the hypothesis implies that
$W\cap\Crit(N_f)=\emptyset$.

For $n\in\N$, we define maps $\eta_n \,:\,W\longrightarrow Y$ as
follows: given $z\in W$, choose a curve $\sigma\in S_{W',W}$ and
$0\leq t\leq 1$ such that $\sigma(t)=z$ and define
$\eta_n(z)=\ext{\sigma}(t+n)$. It is easy to check that this is
well-defined, i.e. independent of the choice of curve and its
parametrization. The maps $\eta_n$ are continuous for all $n \in \N$ with
$\eta_n^{-1}=N_f^{\circ n}$, so all $\eta_n$ are holomorphic.
The sequence $(\eta_n)_{n \in \N}$ clearly forms a normal family,
so there is a locally uniformly convergent subsequence
$(\eta_{n_l})_{l \in \N}$; its limit function $\eta\colon W\to
\ovl Y$ is holomorphic by the theorem of Weierstra{\ss}.

Suppose there is a $y_0 \in \image{\eta}\cap\C$. We have $W\cap
S_0 =\emptyset$, so for all $m \in \N$ and sufficiently big $n$
(depending on $m$), $\image{\eta_n} \cap S_m = \emptyset$.
Furthermore $U=\bigcup_{m \ge 0} S_m$ is open, so $y_0
\in\partial U$ and by Hurwitz' theorem,
$\eta \equiv y_0$ is constant. 
If $y_0\neq\infty$, then choose some
$w\in W'$ and define $z:=N_f(w)$. Then
\begin{equation*}
N_f(y_0)=N_f(\eta(w))=\lim_{l\rightarrow\infty}
N_f(\eta_{n_l}(w))=\lim_{l \rightarrow \infty}
\eta_{n_l}(z)=\eta(z)=y_0 \,\,.
\end{equation*}
Thus $y_0$ is a fixed point of $N_f$ on $\partial U\cap\C$, while
the only fixed points of $N_f$ in $\C$ are the zeros of $f$, and
these are not on the boundary of $U$. This contradiction shows
that $y_0=\infty$. 

Suppose there is a single curve $\sigma\in S_{W,W'}$ for which
$\lim_{t\to\infty}\ext{\sigma}(t)=\infty$ is false. Then there is
a subsequence $(\eta_{n_l})$ for which the limit function
$\eta\colon W\to\ovl Y$ could not be identically equal to
$\infty$, and this is a contradiction which proves the lemma.
\end{proof}

\remark
A different way to prove this lemma is to use contraction
properties of the hyperbolic metric in $U\sm \ovl{P_{U}}$.


\subsection{Immediate Basins are Simply Connected and Unbounded}

Now it is time for our first main result.

\label{Immediate basins:Simple connectivity}
\begin{Theorem}[Immediate basins]
\label{ThmImmediateBasins} \lineclear
Let $f:\C \rightarrow \C$ be a nonlinear entire map, $N_f:=id-\frac{f}{f'}$ its
Newton map and $\xi$ a root of $f$. Then the immediate basin U of $\xi$ is
simply connected and unbounded.
\label{simple_connectivity_immediate_basin}
\end{Theorem}

\begin{proof}
Choose an open disk $S_0$ around $\xi$ with
$N_f(\overline{S_0}) \subset S_0$, $\partial S_0 \cap
P_{U}=\emptyset$ and define $S_{k+1}$ as the component of
$N_f^{-1}(S_k)$ containing
$S_0$. Note that $N_f$ has no postcritical points on $\partial
S_0$, so there are no critical points on any $\partial S_k$ and
$N_f$ is locally biholomorphic in a neighborhood of $\partial
S_k$.

By Lemma \ref{BasinIsCircleInverseImageIB}, 
$U=\bigcup_{n\ge 0} S_n$. If $U$ is not simply connected, then there is
a minimal $M$ such that ${S_M}$ is not full. Choose a bounded
component $B_0$ of $\C \sm S_M$; then there is a bounded component
$B$ of $U\sm S_M$ as well. If $U$ is bounded, let $B:=U\sm S_M$ for an
arbitrary $M$. In both cases, $B$ is a bounded component of $U\sm S_M$,
and this will lead to a contradiction (compare Figure~\ref{B}).

\begin{figure}[h]
\begin{center}
\begin{minipage}[t]{55mm}
\begin{center}\input{B.pstex_t} \end{center}
\caption{$U=\bigcup_{n\ge 0} S_n$ and $B$}
\label{B}
\end{minipage}
\hspace{.5cm}
\begin{minipage}[t]{55mm}
\begin{center}\input{omega.pstex_t}\end{center}
\caption{How to choose $\gamma^0$.}
\label{Omega}
\end{minipage}
\end{center}
\end{figure}

Define 
$P_{B}:={\bigcup_{m\geq 0}N_f^{\circ m}(\Crit(N_f)\cap B)}$; 
since $B$ is bounded, $\Crit(N_f)\cap B$ is finite, and
the only accumulation point of $P_{B}$ is $\xi$. Choose 
$w \in \(S_{M+1} \sm \ovl{S_M}\) \cap B \sm P_B$ (compare
Figure~\ref{Omega}) and set
$z:=N_f(w)\in S_{M} \sm \ovl{S_{M-1}}$. 
There is an injective path $\gamma^0 : [0,1]\rightarrow \(S_{M+1}
\sm \ovl{S_{M-1}}\)\sm P_B$ with $\gamma^0(0)=z$ and
$\gamma^0(1)=w$.

We want to show that $\ext{\gamma^0}$ converges to $\infty$ within
$B$, which would be a contradiction. Since $\gamma^0([0,1])\cap
P_B =\emptyset$, there is a maximal curve
$\gamma:=\ext{\gamma^0}:[0,M_\gamma)\rightarrow \C$ with:
\begin{eqnarray*}
\forall t\in [1,M_\gamma):\quad N_f(\gamma(t))&=&\gamma(t-1)\\
\gamma|_{[0,1]}&\equiv& \gamma^0 \,\,.
\end{eqnarray*}
By induction, it follows that $\gamma([n,n+1])\subset S_{M+n+1} \sm
\ovl{S_{M+n-1}}$ for every $n\in \N$, and in particular we have
$\gamma([1,M_\gamma)) \subset U \setminus \ovl{S_{M}}$. In fact, even
$\gamma([1,M_\gamma))\subset B$ (because $B$ is the
component of $U \setminus S_M$ containing $\gamma(1)=w$).
Since $B$ is bounded, Lemma~\ref{LemExtCurvePossibilities} implies that
$M_\gamma=\infty$.

Choose an open, bounded and simply connected
neighborhood $W \subset S_{M+1} \setminus \ovl{S_{M-1}}$
of $\image{\gamma^0}$ disjoint from $P_B$. This
can be done because $\image{\gamma^0}$ and
$\ovl P_B=P_B\cup\{\xi\}$ are compact and disjoint.

\begin{figure}[h]
{\center \input{wwstrich.pstex_t} \\ } 
\caption[Starting set for curves (Proof of Theorem \ref{simple_connectivity_immediate_basin})]{$W$ and $W'$.}
\label{WWStrich}
\end{figure}

Let $W'\subset W$ be a simply connected neighborhood of $w$ with
$N_f(W')\subset W$
(compare Figure~\ref{WWStrich}). 
By Lemma~\ref{LemExtCurvePossibilities}, every curve
\[
    \sigma \in S:=\{\sigma :[0,1]\rightarrow W \textrm{
    continuous, } \sigma(1) \in W',\, \sigma(0)=N_f(\sigma(1))\}
\]
has an extension $\ext{\sigma}: [0,\infty)\rightarrow B \cup W$.
By Lemma~\ref{BoundedBackwardsImagesOfCurvesConverge}, every
curve $\sigma\in S$ satisfies
$\lim_{t\to\infty}\ext{\sigma}(t)=\infty$, so $B$ is unbounded: a
contradiction. 
\end{proof}

\remark
In many cases, it even follows that $\infty$ is accessible within
$U$. The case in which we cannot prove this is if $U$ contains
infinitely many critical points of $N_f$ such that $P_U$
is dense in $U$.

\section{Virtual Immediate Basins}
\subsection{A Motivating Example}
\label{Virtual basins:Motivation} 

\begin{figure}[h]
{ \center \includegraphics[scale=0.3]{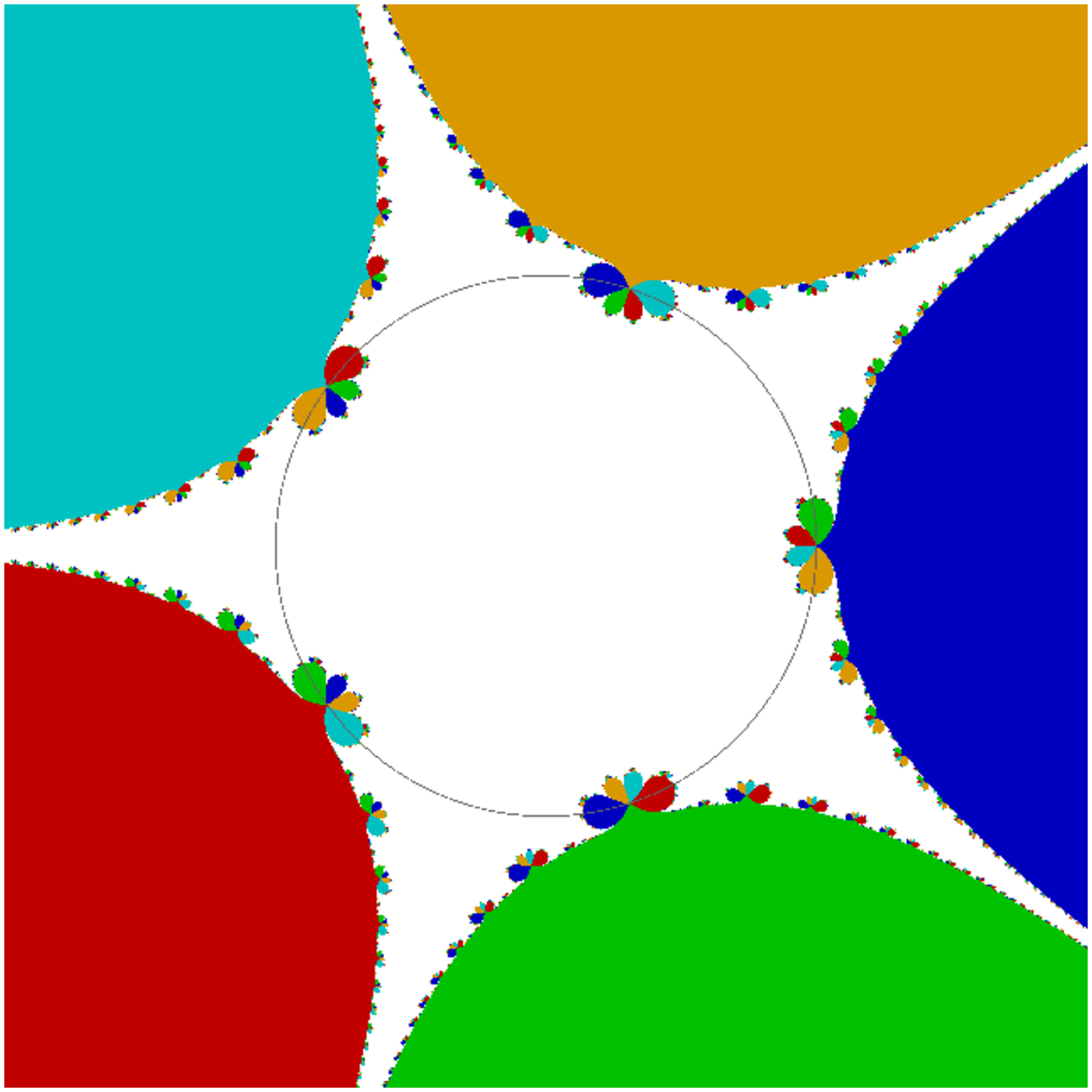}
\includegraphics[scale=0.3]{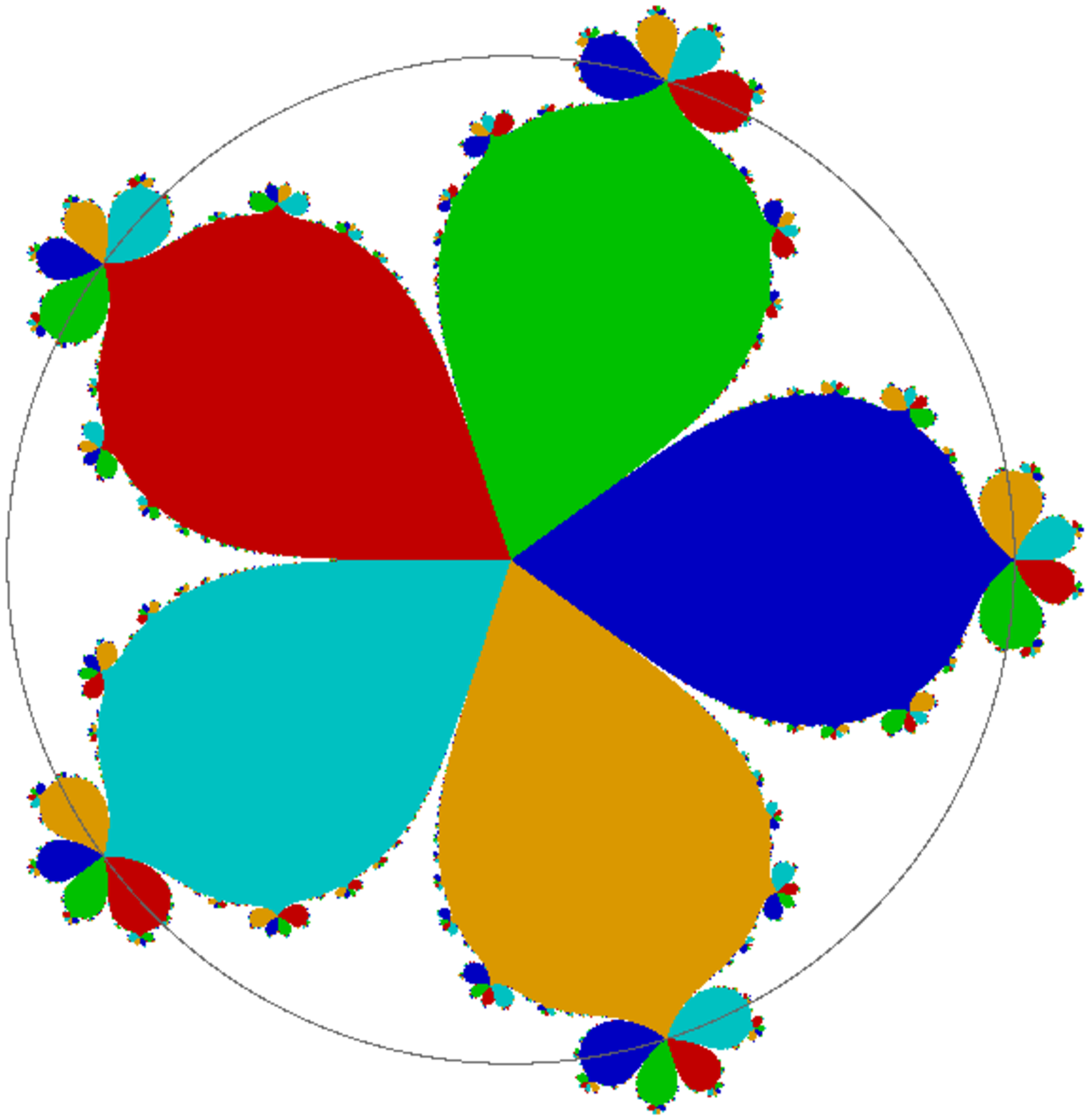}
\\}
\caption[\emph{Left:}
Newton for $f(z)=ze^{-\frac{1}{5}z^5 }$]{Dynamics of the
Newton for $f(z)=ze^{-\frac{1}{5}z^5 }$. The immediate
basin of the root $0$ is white, the other colors correspond to
virtual immediate basins and their backward images. The unit
circle $\S^1$ is marked in grey.  \emph{Right:} The same
situation in $\zeta=\iota(z)=1/z$ coordinates: $\iota \circ N_f
\circ \iota^{-1}(\zeta)=\zeta-\zeta^{n+1}$ is a polynomial.}
\label{figureNewton5}
\end{figure}

The dynamics of Newton's map for transcendental entire functions has a
class of Fatou components which we want to call \emph{virtual immediate
basins}. Let  $f(z):=z \exp({-\frac{1}{n} z^n})$; 
its Newton map $N_f(z)=z \left(1-\frac {1}{1-z^n}\right)$
s a rational function. The involution
$\iota : \hat{\C} \to \hat{\C}, z \mapsto \frac{1}{z}$ conjugates
$N_f$ to the polynomial $\iota \circ N_f \circ
\iota(\zeta)=\zeta-\zeta^{n+1}$. In this case, the Leau-Fatou
``Flower Theorem'' shows that there are
exactly $n$ attracting and $n$ repelling petals at $\zeta=0$, so
$N_f$ has exactly $n$ unbounded Fatou components with convergence
to $z=\infty$; moreover, the immediate basin of the root $0$ has
exactly $n$ accesses to $\infty$ (these accesses are called
``channels to $\infty$'' \cite{HSS}); compare
Figure~\ref{figureNewton5}. We call the attracting petals at
infinity \emph{virtual immediate basins}: their dynamics is
similar as if there was a root at
$\infty$ in each of these $n$ directions. Note that the $n$
channels of the root $0$ separate all these $n$ virtual basins. 

For Newton's method of polynomials, it is known that any pair of
channels to $\infty$ of the same root must enclose a different root of
the polynomial. As this example shows, an analogous statement for
transcendental entire functions would be false if virtual immediate
basins were not taken into account. An investigation of the
combinatorial possibilities between channels to $\infty$ and immediate
basins (including virtual basins) can be found in \cite{JohannesNewton}.

\subsection{An Exhaustion of Virtual Immediate Basins}
\label{Virtual basins:An exhaustion}

In order to define virtual immediate basins, we need the following
definition.
\begin{Definition}[Absorbing set]
  \label{AbsorbingSet}
  If $U$ is an $N_f$-invariant domain in $\C$, then an open set
$A\subset\C$ is called \emph{absorbing set} (of $U$) if
the following conditions hold:
  \begin{enumerate}
  \item $A$ is simply connected and $\overline{A}$ is full;
  \item $\ovl{N_f(A)}\subset A$;
  \item for every $z\in U$ there is a $k$ such that
$N_f^{\circ k}(z)\in A$.
  \end{enumerate}
\end{Definition}

\begin{Definition}[Virtual immediate basin]
  \label{DefVirtualImmediateBasin}
  A domain $U\subset\C$ is called a \emph{virtual immediate
basin} if it is maximal with respect to the
  following properties:
  \begin{enumerate}
  \item $\lim\limits_{n\rightarrow\infty} N_f^{\circ n}(z)=\infty$ for all $z\in U$;
  \item there is an absorbing set $A\subset U$.
  \end{enumerate}
\end{Definition}

Virtual immediate basins occur only for Newton maps of
transcendental entire functions. Obviously, every virtual
immediate basin is unbounded with $\infty$ as accessible boundary
point, and the same is true for every absorbing set of a virtual
immediate basin. If $N_f$ is rational (as in our example above),
then a virtual immediate basin is a Leau-Fatou petal; otherwise,
it is (contained in) a ``Baker domain'' (a domain at
$\infty$); see below. 

Let $U$ be a virtual immediate basin. By a small modification of the
absorbing set, we find an absorbing set $S_0$ of $U$ such that $\partial
S_0$ contains no postcritical points. 
Similarly as for immediate basins, define $S_{k+1}$ to be the connected
component of $N^{-1}(S_k)$ containing $S_0$, for all $k\geq 0$. As
before, we have the following:

\begin{Lemma}
  \label{BasinIsCircleInverseImageVB}
  $U$ is open and $U=\bigcup_{k\in \N} S_k$. If $U$ is not full, then
one of the $S_k$ is not full.
\qed
\end{Lemma}

\subsection{Simple Connectivity}
\label{Virtual basins:Simple connectivity}

In a number of ways, virtual immediate basins have similar
properties as immediate basins of roots; here is one such result.
\begin{Theorem}
\label{SimpleConnectivityVirtualImmediateBasins}
Virtual immediate basins are simply connected.
\end{Theorem}
\begin{proof}
Let $U$ be a virtual immediate basin with absorbing set $S_0$.
By Lemma \ref{BasinIsCircleInverseImageVB} there is an exhaustion
$U=\bigcup_{m\in\N}S_m$ of $U$. If $U$ is not full, then there is a
minimal $M$ such that $S_M$ is not full.

There is a bounded connected component $B$
of $\C \sm S_M$. Its boundary $\partial B$ is a compact subset of $S_M$,
so $N_f(\partial B)$ is a compact subset of $S_{M-1}$.

There are no critical points in $\partial S_M \cup \partial
S_{M-1} \cup \ldots \cup \partial S_0$, so $N_f^M$ maps $\partial
B$ onto a component of $\partial S_0$. Since $\partial B$ is a compact
subset of $\C$, so is $\partial S_0$. But $S_0\subset\C$ is unbounded
and simply connected, so it follows that $\partial S_0$ is unbounded
as well, and this is a contradiction.
\end{proof}

We cannot show in general that every virtual
immediate basin equals an entire Fatou component; however, we
have the following.

\begin{Remark}
Every virtual immediate basin $U$ is contained in an invariant
Fatou component $F$. If $F$ is simply connected and $N|_F$ is a
proper map, then $F$ is a virtual immediate basin.
\end{Remark}
This simply follows by using a Riemann map $\varphi: F\rightarrow \D$ to
transport the dynamics of $F$ into the unit disk $\D$. It might
be possible to extend Shishikura's results~\cite{Mitsu} to the
transcendental case, to show the more general result that every
Fatou component of Newton's map for entire functions is simply
connected. In this case, every virtual immediate basin would be an
entire Fatou component provided the Newton dynamics restricted to
this Fatou component was proper.


\end{document}

%% file: B.pstex_t
\begin{picture}(0,0)%
\includegraphics{B.pstex}%
\end{picture}%
\setlength{\unitlength}{1657sp}%
\begingroup\makeatletter\ifx\SetFigFont\undefined%
\gdef\SetFigFont#1#2#3#4#5{%
  \reset@font\fontsize{#1}{#2pt}%
  \fontfamily{#3}\fontseries{#4}\fontshape{#5}%
  \selectfont}%
\fi\endgroup%
\begin{picture}(6788,4524)(425,-4123)
\put(4951,-2266){\makebox(0,0)[lb]{\smash{\SetFigFont{7}{8.4}{\rmdefault}{\mddefault}{\updefault}{\color[rgb]{0,0,0}$p$}%
}}}
\put(4771,-1771){\makebox(0,0)[lb]{\smash{\SetFigFont{7}{8.4}{\rmdefault}{\mddefault}{\updefault}{\color[rgb]{0,0,0}$S_0$}%
}}}
\put(4006,-1636){\makebox(0,0)[lb]{\smash{\SetFigFont{7}{8.4}{\rmdefault}{\mddefault}{\updefault}{\color[rgb]{0,0,0}$S_1$}%
}}}
\put(1441,-1366){\makebox(0,0)[lb]{\smash{\SetFigFont{7}{8.4}{\rmdefault}{\mddefault}{\updefault}{\color[rgb]{0,0,0}$S_M$}%
}}}
\put(3421,-2266){\makebox(0,0)[lb]{\smash{\SetFigFont{7}{8.4}{\rmdefault}{\mddefault}{\updefault}{\color[rgb]{0,0,0}$\dots$}%
}}}
\put(3016, 29){\makebox(0,0)[lb]{\smash{\SetFigFont{7}{8.4}{\rmdefault}{\mddefault}{\updefault}{\color[rgb]{0,0,0}$\C \setminus U$}%
}}}
\put(3286,-4021){\makebox(0,0)[lb]{\smash{\SetFigFont{7}{8.4}{\rmdefault}{\mddefault}{\updefault}{\color[rgb]{0,0,0}$\C \setminus U$}%
}}}
\put(2431,-1231){\makebox(0,0)[lb]{\smash{\SetFigFont{7}{8.4}{\rmdefault}{\mddefault}{\updefault}{\color[rgb]{0,0,0}$S_{M-1}$}%
}}}
\put(1711,-2176){\makebox(0,0)[lb]{\smash{\SetFigFont{7}{8.4}{\rmdefault}{\mddefault}{\updefault}$B_0$}}}
\put(1846,-3391){\makebox(0,0)[lb]{\smash{\SetFigFont{7}{8.4}{\rmdefault}{\mddefault}{\updefault}{\color[rgb]{0,0,0}$B$}%
}}}
\end{picture}

%% file: omega.pstex_t
\begin{picture}(0,0)%
\includegraphics{omega.pstex}%
\end{picture}%
\setlength{\unitlength}{1657sp}%
\begingroup\makeatletter\ifx\SetFigFont\undefined%
\gdef\SetFigFont#1#2#3#4#5{%
  \reset@font\fontsize{#1}{#2pt}%
  \fontfamily{#3}\fontseries{#4}\fontshape{#5}%
  \selectfont}%
\fi\endgroup%
\begin{picture}(6788,4524)(425,-4123)
\put(4951,-2266){\makebox(0,0)[lb]{\smash{\SetFigFont{7}{8.4}{\rmdefault}{\mddefault}{\updefault}{\color[rgb]{0,0,0}$p$}%
}}}
\put(1441,-1366){\makebox(0,0)[lb]{\smash{\SetFigFont{7}{8.4}{\rmdefault}{\mddefault}{\updefault}{\color[rgb]{0,0,0}$S_M$}%
}}}
\put(3736,-3031){\makebox(0,0)[lb]{\smash{\SetFigFont{7}{8.4}{\rmdefault}{\mddefault}{\updefault}{\color[rgb]{0,0,0}$z$}%
}}}
\put(3106,-3256){\makebox(0,0)[lb]{\smash{\SetFigFont{7}{8.4}{\rmdefault}{\mddefault}{\updefault}{\color[rgb]{0,0,0}$\gamma_0$}%
}}}
\put(2566,-1636){\makebox(0,0)[lb]{\smash{\SetFigFont{7}{8.4}{\rmdefault}{\mddefault}{\updefault}{\color[rgb]{0,0,0}$S_{M+1}$}%
}}}
\put(3151, 74){\makebox(0,0)[lb]{\smash{\SetFigFont{7}{8.4}{\rmdefault}{\mddefault}{\updefault}{\color[rgb]{0,0,0}$\C \setminus U$}%
}}}
\put(3196,-3976){\makebox(0,0)[lb]{\smash{\SetFigFont{7}{8.4}{\rmdefault}{\mddefault}{\updefault}{\color[rgb]{0,0,0}$\C \setminus U$}%
}}}
\put(1756,-2131){\makebox(0,0)[lb]{\smash{\SetFigFont{7}{8.4}{\rmdefault}{\mddefault}{\updefault}$B_0$}}}
\put(1486,-3301){\makebox(0,0)[lb]{\smash{\SetFigFont{7}{8.4}{\rmdefault}{\mddefault}{\updefault}{\color[rgb]{0,0,0}$B$}%
}}}
\put(3106,-2851){\makebox(0,0)[lb]{\smash{\SetFigFont{7}{8.4}{\rmdefault}{\mddefault}{\updefault}{\color[rgb]{0,0,0}$w$}%
}}}
\end{picture}

%% file: wwstrich.pstex_t
\begin{picture}(0,0)%
\includegraphics{wwstrich.pstex}%
\end{picture}%
\setlength{\unitlength}{2368sp}%
\begingroup\makeatletter\ifx\SetFigFont\undefined%
\gdef\SetFigFont#1#2#3#4#5{%
  \reset@font\fontsize{#1}{#2pt}%
  \fontfamily{#3}\fontseries{#4}\fontshape{#5}%
  \selectfont}%
\fi\endgroup%
\begin{picture}(6398,3174)(13,-2319)
\put(2851,-211){\makebox(0,0)[lb]{\smash{\SetFigFont{8}{9.6}{\rmdefault}{\mddefault}{\updefault}{\color[rgb]{0,0,0}$N_f$}%
}}}
\put(676,-1036){\makebox(0,0)[lb]{\smash{\SetFigFont{8}{9.6}{\rmdefault}{\mddefault}{\updefault}{\color[rgb]{0,0,0}$W'$}%
}}}
\put(2401,-2236){\makebox(0,0)[lb]{\smash{\SetFigFont{8}{9.6}{\rmdefault}{\mddefault}{\updefault}{\color[rgb]{0,0,0}$S_M$}%
}}}
\put(5476,-286){\makebox(0,0)[lb]{\smash{\SetFigFont{8}{9.6}{\rmdefault}{\mddefault}{\updefault}{\color[rgb]{0,0,0}$z$}%
}}}
\put(751,-211){\makebox(0,0)[lb]{\smash{\SetFigFont{8}{9.6}{\rmdefault}{\mddefault}{\updefault}{\color[rgb]{0,0,0}$w$}%
}}}
\put(2926,-1036){\makebox(0,0)[lb]{\smash{\SetFigFont{8}{9.6}{\rmdefault}{\mddefault}{\updefault}{\color[rgb]{0,0,0}$\gamma^0$}%
}}}
\put(3151,-1636){\makebox(0,0)[lb]{\smash{\SetFigFont{8}{9.6}{\rmdefault}{\mddefault}{\updefault}{\color[rgb]{0,0,0}$W$}%
}}}
\end{picture}